\documentclass[10pt]{article} 
\usepackage{amsfonts} 
 
\usepackage{amscd}

\newtheorem{theorem}{Theorem}[section]  
  
\newtheorem{lemma}[theorem]{Lemma}  
\newtheorem{proposition}[theorem]{Proposition}  
  
\newtheorem{definition}{Definition}[section]

\begin{document} 
\title{Lie-admissible algebras and Operads}
\author{Michel Goze . Elisabeth Remm}
\maketitle 
\begin{center}
Universit\'{e} de Haute Alsace, F.S.T.\\ 
4, rue des Fr\`{e}res Lumi\`{e}re - 68093 MULHOUSE - France 
\end{center}
\begin{abstract}
A Lie-admissible algebra gives a Lie algebra by anticommutativity. In this work we describe remarkable types of Lie-admissible
 algebras such as Vinberg algebras, pre-Lie algebras
 or Lie algebras. We compute the corresponding binary quadratic operads and study their duality.
 Considering Lie algebras as Lie-admissible algebras, we can define for each
Lie algebra a cohomology with values in an Lie-admissible module. This permits to
study some deformations of Lie algebras, in particular classes of Lie-admissible
 algebras such as Vinberg algebras or pre-Lie algebras.
\end{abstract}

\section{Lie-admissible algebras}

\subsection{Definition}

Let $\mathcal{A}=(A,\mu )$ be a finite dimensional algebra on a commutative
field $\Bbb{K}$ of characteristic zero. In this notation, $\mu $ is the law
of $\mathcal{A}$, that is a linear mapping 
\[
\mu :A\otimes A\rightarrow A 
\]
on the vector space $A.$

We denote by $a_{\mu }:A^{\otimes 3}\rightarrow A$
the associator of the law  $\mu $ :  
\[
a_{\mu }\left( X_1,X_2,X_3\right) =\mu \left( \mu \left( X_1,X_2\right)
,X_3\right) -\mu \left( X_1,\mu \left( X_2,X_3\right) \right).
\]
Let $\sum_{n}$ be the symmetric group of degree $n$. For every $\sigma \in
\sum_{3}$, we put
\[
\sigma \left( X_{1},X_{2},X_{3}\right) = \left( X_{\sigma ^{-1}\left( 1\right) },X_{\sigma ^{-1}\left( 2\right)
},X_{\sigma ^{-1}\left( 3\right) }\right).
\]
\begin{definition}
The algebra $\mathcal{A}=(A,\mu )$ is called Lie-admissible if the
law $\mu $ satisfies 
\[
\sum_{\sigma \in \sum_{3}}\left( -1\right) ^{\varepsilon \left( \sigma
\right) }a_{\mu } \circ \sigma=0. \qquad \qquad (\ast )
\]
\end{definition}
This definition ([A]) is equivalent to say that the mapping 
$[,]_{\mu}:A\otimes A\rightarrow A$ defined by $[X,Y]=\mu (X,Y)-\mu (Y,X)$ is a Lie
bracket. We will denote by $\mathcal{A}_{L}$ the corresponding Lie algebra.
\subsection{Examples}

\noindent 1. Every associative algebra (not necessarily unitary) is a Lie-admissible
 algebra.

\noindent 2. An algebra $\mathcal{A}=(A,\mu )$ is a Vinberg algebra (also
called left
 symmetric algebra) if its law satisfies 
\[
\mu (X,\mu (Y,Z))-\mu (Y,\mu (X,Z))=\mu (\mu (X,Y),Z)-\mu (\mu (Y,X),Z). 
\]
It is also a Lie-admissible algebra.

\noindent 3. Of course a Lie algebra law is a law of Lie-admissible
algebra, the Jacobi conditions implying $(\ast ).$

\noindent 4. A pre-Lie algebra is defined by a law $\mu $ such that 
\[
\mu (\mu (X,Y),Z)-\mu (X,\mu (Y,Z))=\mu (\mu (X,Z),Y)-\mu (X,\mu (Z,Y)). 
\]
The bracket $[,]_{\mu }$ being a Lie bracket, a pre-Lie algebra is a Lie-admissible algebra.

\noindent \textbf{Remarks.} Every associative algebra is a Vinberg algebra and a pre-Lie algebra. A
Vinberg algebra is a pre-Lie algebra if and only if all the associators are
equal  
\[
a_\mu \circ \sigma=a_\mu \circ \upsilon 
\]
for every $\upsilon $ and $\sigma $ in $\sum_{3}.$

A Lie algebra is associative (respectively Vinberg, respectively pre-Lie
algebra) if and only if it is 2-step nilpotent.

\subsection{Geometrical interpretation}

Let $\frak{g}$ be a real Lie algebra provided with an affine structure. The
associated covariant derivative $\bigtriangledown $ satisfies 
\[
\left\{ 
\begin{array}{l}
\bigtriangledown _{X}Y-\bigtriangledown _{Y}X-[X,Y]=0 \\ 
\bigtriangledown _{X}\bigtriangledown _{Y}-\bigtriangledown
_{Y}\bigtriangledown _{X}=\bigtriangledown _{[X,Y]}
\end{array}
\right. 
\]
Then the product $\mu (X,Y)=\bigtriangledown _{X}Y$ endows the vector space $%
\frak{g}$ with a structure of Vinberg algebra. As we
have $  [X,Y]=\mu (X,Y)-\mu (Y,X),$ every Lie algebra provided with an affine
structure is subordinated to a Lie-admissible algebra which is in this case
a Vinberg algebra. More generally, every Lie algebra is subordinated to a
Lie-admissible algebra.\ In fact it is sufficient to consider Levi Civita
connections (which always exists). As the torsion vanishes, the covariant derivative satisfies 
the first of the previous equations and $\mu(X,Y)=\bigtriangledown _{X}Y$ is a law af Lie-admissible algebra such that 
$[X,Y]=[X,Y]_{\mu }.$

This permits to consider the set of Lie-admissible algebras as the set of
invariant linear connections torsionless on Lie algebras. In fact,
if $\mu $ is a law of admissible algebra, putting $\bigtriangledown
_{X}Y=\mu (X,Y)$ then $\bigtriangledown $ defines a linear connection if
the Bianchi identities are satisfied. Let us note $R(X,Y)$ the curvature
tensor corresponding to $\bigtriangledown .\;$As $\mu $ is an Lie-admissible
law, we have $R(X,Y).Z+R(Y,Z).X+R(Z,X).Y=0$ then the first identity is
realized. For the second we can prove that 
\[
(\bigtriangledown _{X}R)(Y,Z)+(\bigtriangledown
_{Y}R)(Z,X)+(\bigtriangledown _{Z}R)(X,Y)=0. 
\]
Using the relations $(\bigtriangledown _{X}R)(Y,Z)=\bigtriangledown
_{X}(R(Y,Z))-R(\bigtriangledown _{X}Y,Z)-R(Y,\bigtriangledown _{X}Z),$ and $%
\bigtriangledown _{[X,Y]}=\bigtriangledown _{\bigtriangledown
_{X}Y}-\bigtriangledown _{\bigtriangledown _{Y}X}$, we deduce the second
identity.

\subsection{Actions of the symmetric group $\sum_{3}$}

\begin{definition}
Let $G$ be a subgroup of ${\sum}_{3}.$ We say that the algebra law is $G$-associative if 
\[
\sum_{\sigma \in G}(-1)^{\varepsilon (\sigma )}a_{\mu }\circ \sigma
=0. 
\]
\end{definition}
The subgroups of ${\sum}_{3}$ are
well known. We have $G_{1}=\{Id\},G_{2}=\{Id,\tau _{12}\},G_{3}=\{Id,\tau
_{23}\},G_{4}=\{Id,\tau _{13}\},G_{5}=\{Id,(231),(312)\}=A_3$ (the alternating group), $G_{6}=\sum_{3}$
where $\tau _{ij}$ is the transposition between $i$ and $j$ and $(231)$ a
$3$-cycle. We deduce the following type of Lie-admissible algebras :

1. If $\mu $ is $G_{1}$-associative then $\mu $ is an associative law.

2. If $\mu $ is $G_{2}$-associative then $\mu $ is a law of Vinberg algebra (also called left symmetric algebras [N]). 
If $A$ is finite-dimensional, the associated Lie algebra is provided with an affine structure.

3. If $\mu $ is $G_{3}$-associative then $\mu $ is a law of pre-Lie algebra [G], also called right symmetric algebras. 
Let us note that if $x.y$ is pre-Lie law then $x \odot y=y.x$ is a Vinberg law.
 
4. If $\mu $ is $G_{4}$-associative then $\mu $ satisfies 
\[
(X.Y).Z-X.(Y.Z)=(Z.Y).X-Z.(Y.X) 
\]

5. If $\mu $ is $G_{5}$-associative then $\mu $ satisfies the generalized
Jacobi condition : 
\[
(X.Y).Z+(Y.Z).X+(Z.X).Y=X.(Y.Z)+Y.(Z.X)+Z.(X.Y) 
\]
Moreover if the law is antisymmetric, then it is a law of Lie algebra.

6. If $\mu $ is $G_{6}$-associative then $\mu $ is a Lie-admissible law.

\smallskip

\noindent{\bf{Example. }}
Let us consider the vector space $\mathcal{C}^{\infty }(\Bbb{R}$,$\Bbb{R)}$
of real infinitely derivable functions with values in $\Bbb{R}$.
We can endow this space with the following algebra structures :

1. $\mu _{1}(f,g)=f.g$ and $(\mathcal{C}^{\infty }(\Bbb{R}$,$\Bbb{R}$),$\mu
_{1})$ is associative (type 1)

2. $\mu _{2}(f,g)=f.g^{\prime }$ and $(\mathcal{C}^{\infty }(\Bbb{R}$,$\Bbb{R%
}$),$\mu _{2})$ is a Vinberg algebra (type 2)

3. $\mu _{3}(f,g)=f^{\prime }.g$ and $(\mathcal{C}^{\infty }(\Bbb{R}$,$\Bbb{R%
}$),$\mu _{3})$ is a pre-Lie algebra (type 3)

4. $\mu _{4}(f,g)=f^{\prime }.g+f.g^{\prime }$ and $(\mathcal{C}^{\infty }(%
\Bbb{R}$,$\Bbb{R}$),$\mu _{4})$ is invariant by $G_{4}$ (type 4).

Suppose that $\mu $ is invariant by $G_{5}.$ If this law is also
antisymmetric then it is a Lie algebra law. This shows that the definition
of $G_{5}$-invariance gives another generalization of ''non-commutative''
Lie algebra than the notion of Leibniz algebra. 

\subsection{$G$-cogebras}
In this section we  introduce the notion of cogebra dualizing the $G_i$-associative algebras. This leads to present the definition in the 
category $Vect_{\Bbb{K}}$ of $\Bbb{K}$-vector spaces.
\begin{definition}
A $G_i$-associative algebra is a pair $(A,\mu)$ where $A$ is a vector space and 
$\mu :A \otimes A : \longrightarrow A $ a linear mapping satisfying the 
following axiom ($G_i$-ass):  

\noindent The square
$$
\begin{CD}
A \otimes A  \otimes A @>{(\mu \otimes Id})_{G_i}>> A \otimes A \\
@V{(Id \otimes \mu)}_{G_i}VV  @V{\mu}VV \\
A \otimes A  @ >{\mu}>> A
\end{CD}
$$
commutes, where ${(Id \otimes \mu)}_{G_i}$ is the linear mapping defined by :
$$
{(Id \otimes \mu)}_{G_i}=\sum_{\sigma \in G_i} (-1)^{\epsilon (\sigma )} (Id \otimes \mu) \circ \sigma.
$$
\end{definition}
Let $\Delta$ be a comultiplication on a vector space $C$:
$$\Delta : C \longrightarrow C \otimes C.$$
 We define the bilinear mapping $G_i \circ (\Delta \otimes Id)$ by 
$$G_i \circ (\Delta \otimes Id)=\sum_{\sigma \in G_i} (-1)^{\epsilon (\sigma )} \sigma \circ (\Delta \otimes Id).$$
\begin{definition}
A $G_i$-cogebra is a vector space $C$ provided with a comultiplication $\Delta : C \longrightarrow C \otimes C$ such that the following diagram is commutative :
$$
\begin{CD}
C  @>{\Delta}>> C \otimes C \\
@V{\Delta}VV  @VV{G_i \circ (Id \otimes \Delta)}V \\
C \otimes C  @ >{G_i \circ (\Delta \otimes Id)}>> C \otimes C \otimes C.
\end{CD}
$$
\end{definition}
The next results relate $G_i$-algebras and $G_i$-cogebras.
\begin{proposition}
The dual space of a $G_i$-cogebra is a $G_j$-associative algebra .
\end{proposition}
{\it Proof}. Let $(C,\Delta)$ a $G_i$-cogebra. Let us consider the maps $\lambda : C^{*} \otimes C^{*} \rightarrow (C \otimes C)^{*}$ given 
by $\lambda(f \otimes g)(v \otimes u)=f(u) \otimes g(v)$ and $\tau : C^{*} \otimes C^{*} \rightarrow C^{*} \otimes C^{*}$ defined by 
$\tau (u \otimes v)=v \otimes u$. Let us consider the law $\mu=\Delta ^{*} \circ \lambda \circ \tau$. It provides the dual space $C^{*}$ 
with a $G_j$-associative algebra. In fact we have
$$\mu (f_1 \otimes f_2)(x)=\Delta ^{*}\circ \lambda (f_2 \otimes f_1)(x)=(f_2 \otimes f_1)(\Delta (x))$$
for all $f_1 , f_2 \in C^{*}, x \in C.$ We denote $X \otimes Y = \Delta (x), X_1 \otimes X_2 = \Delta (X)$ and $Y_1 \otimes Y_2 = \Delta (Y)$ (we use the sigma notation and we forget the sigma). This gives 
$$ \mu (f_1 \otimes f_2)(x)=f_1 (X)f_2 (Y).$$
The associator of $\mu$ is written :
$$\mu (\mu (f_1 \otimes f_2) \otimes f_3)(x)-\mu (f_1 \otimes \mu ( f_2 \otimes f_3))(x)=(f_1 \otimes f_2 \otimes f_3)(X_1 \otimes X_2 \otimes Y -X \otimes Y_1 \otimes Y_2)$$
and the property of the comultiplication $\Delta$ is equivalent to
$$ \sum_{\sigma \in G_i} (-1)^{\epsilon (\sigma )} \sigma \circ (\Delta \otimes Id)(\Delta (x))=\sum_{\sigma \in G_i} (-1)^{\epsilon (\sigma )} \sigma \circ (Id \otimes \Delta )(\Delta (x))$$ that is
$$\sum_{\sigma \in G_i} (-1)^{\epsilon (\sigma )} \sigma  (\Delta (X) \otimes Y)=\sum_{\sigma \in G_i} (-1)^{\epsilon (\sigma )} \sigma  (X \otimes \Delta (Y))$$
or
$$\sum_{\sigma \in G_i} (-1)^{\epsilon (\sigma )} \sigma (X_1 \otimes X_2 \otimes Y - X \otimes Y_1 \otimes Y_2)=0. $$
This relation proves that the associator of $\mu$ satisfies the same relation. 
\begin{proposition}
The dual vector space of a finite dimensional $G_i$-associative algebra has a $G_i$-cogebra structure.
\end{proposition}
{\it Proof.} Let $A$ a finite dimensional $G_i$-associative algebra and let $\{{e_i},{i=1,...,n}\}$ be a basis of $A$. If 
$\{f_i\}$ is the dual basis then $\{f_i \otimes f_j\}$ is a basis of $A^* \otimes A^* $. Let us put
$$\Delta (f)=\sum _{i,j} f(\mu (e_i \otimes e_j)) f_i \otimes f_j.$$
In particular
$$\Delta (f_k) = \sum _{i,j} C_{ij}^{k}f_i \otimes f_j$$
where $C_{ij}$ are the structure constants of $\mu$ related to the basis $\{{e_i}\}$. Then $\Delta$ is the comultiplication of a $G_i$-cogebra.

\section{Lie-Admissible operads}

\subsection{Binary quadratic operads}

Let $\Bbb{K} [ \sum_n ]$ be the $\Bbb{K}$-algebra of the symmetric group $\sum_n$.
An operad $\mathcal{P}$ is defined by a sequence of $\Bbb{K}$-vector spaces $\mathcal{P}%
(n)$, $n\geq 1$ such that $\mathcal{P}(n)$ is a module over $\Bbb{K}\left[
\sum_{n}\right] $ and
with composition maps 
\[
\circ _{i}:\mathcal{P}(n)\otimes \mathcal{P}(m)\rightarrow \mathcal{P}%
(n+m-1)\qquad i=1,...,n 
\]
satisfying some ''associative'' properties, the May Axioms [M].

Any $\Bbb{K}\left[ \sum_{n}\right] $-module $E$ generates a free operad noted $\mathcal{F(}E)$ ([G])
satisfying $\mathcal{F(}E)(1)=\Bbb{K}$, $\mathcal{F(}E)(2)=E$. In particular if $E=\Bbb{K}\left[ \sum_{2}\right]$, the free module $\mathcal{F(}E)(n)$ admits as a basis the
"parenthized products" of $n$ variables indexed by $\{1,2,...,n\}$. For instantce a basis of 
$\mathcal{F(}E)(2)$ is given by $(x_{1}.x_{2})$ and $(x_{2}.x_{1})$, and a basis of 
$\mathcal{F(}E)(3)$ is given by 
$$\left\{ \left( (x_{i}.x_{j}).x_{k}\right)
,\left( x_{i}.(x_{j}.x_{k}).\right) ,i\neq j\neq k\neq i,i,j,k\in \left\{
1,2,3\right\} \right\} .$$
Let $E$ be a $\Bbb{K}\left[ \sum_{2}\right]$-module and $R$ a $\Bbb{K}\left[ \sum_{3}\right]$-submodule of 
$\mathcal{F(}E)(3)$. We denote $\mathcal{R}$ the ideal generated by $R$, that is the intersection of all the ideals 
$\mathcal{I}$ of $\mathcal{F}(E)$ such that $\mathcal{I}(1)=0,\mathcal{I}(2)=0$ and 
$\mathcal{I}(3)=R$.  

We call binary quadratic operad generated by $E$ and $R$ the operad $\mathcal{P}(\Bbb{K},E,R)$, also
denoted  $\mathcal{F(}E)/\mathcal{R}$ and defined by 
\[
\mathcal{P}(\Bbb{K},E,R)(n)=
(\mathcal{F(}E)/\mathcal{R)}(n)=\frac{\mathcal{F(}E)(n)}{\mathcal{R(}n%
\mathcal{)}} 
\]

\noindent Thus an operad $\mathcal{P}$ is binary quadratic operad if and only if there exists a $%
\Bbb{K}\left[ \sum_{2}\right] $-submodule $E$ and $R$ a $\Bbb{K}\left[
\sum_{3}\right] $-submodule of $\mathcal{F(}E)(3)$ such that $%
\mathcal{P\simeq F(}E)/\mathcal{R}$.

\smallskip

\noindent {\bf Example}. The associative operad $\mathcal{A}ss,$ the Lie operad $\mathcal{L}%
ie$, the Leibniz operad $\mathcal{L}eib$ are binary quadratic operads ([G.K],[L]).

\subsection{Lie-Admissible operads}

Let $\mathcal{F(}E)$ be the free operad generated by $E=\Bbb{K}\left[ \sum_{2}\right]$. 
Consider $R$ the $\Bbb{K}\left[ \sum_{3}\right] $-submodule generated by the
vector 
\begin{eqnarray*}
u
&=&x_{1}.(x_{2}.x_{3})+x_{2}.(x_{3}.x_{1})+x_{3}.(x_{1}.x_{2})-x_{2}.(x_{1}.x_{3})-x_{3}.(x_{2}.x_{1}) \\
&&-x_{1}.(x_{3}.x_{2})
-(x_{1}.x_{2}).x_{3}-(x_{2}.x_{3}).x_{1}-(x_{3}.x_{1}).x_{2}+(x_{2}.x_{1}).x_3 \\
&& +(x_{3}.x_{2}).x_{1}+(x_{1}.x_{3}).x_{2}
\end{eqnarray*}

\noindent The Lie-Admissible operad, noted $\mathcal{L}ieAdm$ is the binary quadratic operad
defined by 
\[
\mathcal{L}ieAdm=\mathcal{F(}E)/\mathcal{R} 
\]
As we have distinguished 6 types of Lie-admissible algebras, we can
determine each corresponding operad. We obtain the
 following binary quadratic
operads :
 
1. $\mathcal{A}ss$ corresponding to the associative algebras.

2. $\mathcal{V}inb$ corresponding to the Vinberg algebras. Here the $\Bbb{K}%
\left[ \sum_{3}\right] $-submodule is generated by the vectors 
\[
x_{1}.(x_{2}.x_{3})-(x_{1}.x_{2}).x_{3}-x_{2}.(x_{1}.x_{3})+(x_{2}.x_{1}).x_{3}. 
\]

3. $\mathcal{P}reLie$ corresponding to the pre-Lie algebras ($G_{3}$%
-associative). The vectors which generate the ideal $R$ are : 
\[
x_{1}.(x_{2}.x_{3})-(x_{1}.x_{2}).x_{3}-x_{3}.(x_{2}.x_{1})+(x_{3}.x_{2}).x_{1}. 
\]

4. $G_{4}-\mathcal{A}ss$ corresponding to the $G_{4}-$associative algebras. $%
R$ is generated by 
\[
x_{1}.(x_{2}.x_{3})-(x_{1}.x_{2}).x_{3}-x_{3}.(x_{2}.x_{1})+(x_{3}.x_{2}).x_{1}. 
\]

5. $G_{5}-\mathcal{A}ss$ corresponding to the $G_{5}-$associative algebras. $%
R$ is generated by 
\[
x_{1}.(x_{2}.x_{3})-(x_{1}.x_{2}).x_{3}+x_{2}.(x_{3}.x_{1})-(x_{2}.x_{3}).x_{1}+x_{3}.(x_{1}.x_{2})-(x_{3}.x_{1}).x_{2}. 
\]

6. $\mathcal{L}ieAdm.$

\subsection{The dual operad $\mathcal{L}ieAdm^{!}$}

Let us consider the binary quadratic operad 
$\mathcal{P}(\Bbb{K},E,R)$, the dual binary quadratic
operad is 
$$\mathcal{P}^{!}=\mathcal{P}(\Bbb{K}^{op},E^{\vee },R^{\perp })$$
where $E^{\vee }$ is the dual of $E$ tensorised by the signature of $\sum_n$ and $R^{\perp }$ the orthogonal
complement to $R$ in $\mathcal{F(}E^{\vee })(3)=\mathcal{F(}E)(3)^{\vee }$.

\noindent Before studying the dual operads of each Lie-Admissible operad defined above,
let us introduce some classes of associative algebras.

\begin{definition}
An associative algebra $A$ is called 3-order abelian if we have 
\[
X_{1}X_{2}X_{3}=X_{\sigma (1)}X_{\sigma (2)}X_{\sigma (3)}
\]
for all $\sigma \in \sum_{3}$ and for all $X_{i}\in A.$
\end{definition}

The unitary 3-order abelian algebras are the commutative algebras. But there
exists non commutative 3-order abelian algebras. Let us consider, for
example, the five dimensional associative algebra defined by 
\[
e_{1}^{2}=e_{2},\quad
e_{1}^{3}=e_{3}, \quad
e_{1}e_{4}=e_{5}, \quad e_{4}e_{1}=e_{3}+e_{5}, \quad e_{4}^{2}=e_{3}. 
 \]
If $A$ is a 3-order abelian algebra, then the subalgebra $\mathcal{D}(A)$
generated by the product $xy$ is abelian.\ Then $A$ is an extension
\[
0\rightarrow V\rightarrow A\rightarrow A_{1}\rightarrow 0 
\]
where $A_{1}$ is abelian and $V$ satisfying $vx=xv$ for all $x\in A_{1}$ and 
$v\in V.$ In this case, the corresponding Lie algebra is 2-step nilpotent.
In fact, as we have $abc=bac$, then $[a,b]c=0.$ We have also $c[a,b]=0$ thus 
$[[a,b],c]=0.$ Return to the geometrical interpretation. Such a Lie algebra
is provided with an affine connection which have neither curvature nor torsion.
Moreover the operator of connection satisfies 
\[
\nabla _{X}=0 
\]
for all $X$ in the center and 
\[
\nabla _{X}\nabla _{Y}=\nabla _{Y}\nabla _{X} 
\]
for any generators $X$ and $Y$.

\noindent Let us consider now the scalar product on $\mathcal{F(}E)(3)$ defined by 
\begin{eqnarray*}
&<i(jk),i(jk)>=sgn(
\begin{tabular}{lll}
$1$ & $2$ & $3$ \\ 
$i$ & $j$ & $k$%
\end{tabular}
) \\
&<(ij)k,(ij)k>=-sgn(
\begin{tabular}{lll}
$1$ & $2$ & $3$ \\ 
$i$ & $j$ & $k$%
\end{tabular}
)
\end{eqnarray*}
Let $R$ be the $\Bbb{K}\left[ \sum_{3}\right] $-submodule which determines the
$\mathcal{L}ieAdm$ operad. The annihilator $R^{\perp }$ respect to this scalar
product is of dimension $11$. Let $R^{\prime }$ be the $\Bbb{K}\left[
\sum_{3}\right] $-submodule of $\mathcal{F(}E)(3)$ generated by the
relations 
\begin{eqnarray*}
(x_{\sigma (1)}x_{\sigma (2)})x_{\sigma (3)}-x_{\sigma (1)}(x_{\sigma
(2)}x_{\sigma (3)}),\\
(x_{\sigma (1)}x_{\sigma (2)})x_{\sigma (3)}-(x_{\sigma
(1)}x_{\sigma (3)})x_{\sigma (2)}, \\
(x_{\sigma (1)}x_{\sigma (2)})x_{\sigma (3)}-(x_{\sigma (2)}x_{\sigma
(1)})x_{\sigma (3)}.
\end{eqnarray*}
Then $\dim R^{\prime }=11$ and $<u,v>=0$ for all $v\in R^{\prime }$ where $u$
is the vector which generates $R.$ This implies $R^{\prime }\simeq R^{\perp }
$ and ($\mathcal{F(}E)/\mathcal{R})^{!}$ is by definition to binary the quadratic operad 
$\mathcal{F(}E)/\mathcal{R}^{\perp }.$

\begin{proposition}
The dual operad of $\mathcal{L}ieAdm$ is the binary quadratic operad whose corresponding algebras 
are associative and satisfying 
\[
abc=acb=bac
\]
that is there are 3-order abelian.
\end{proposition}

Return to the different classes of Lie-admissible algebras. One knows that $%
\mathcal{A}ss^{! }=\mathcal{A}ss$ ([G.K]) and $\mathcal{P}reLie^{! }=%
\mathcal{P}erm$ ([Ch]). Recall that this last corresponds to the associative
product with the identities :
\[
abc=acb
\]

\begin{proposition}
The dual operads of $\mathcal{V}inb,$ $G_{4}-\mathcal{A}ss,$ $G_{5}-\mathcal{%
A}ss$ are the quadratic operads with corresponding algebras the associative algebras
satisfying respectively :

- for $\mathcal{V}inb^{! }$ : $abc=bac$

- for $G_{4}-\mathcal{A}ss^{!}$ : $abc=cba$

- for $G_{5}-\mathcal{A}ss^{! }$ : $abc=bca=cab.$
\end{proposition}

\noindent {\it Sketch of proof.} $R_{2}^{\perp }$ is the 
$\Bbb{K}\left[ \sum_{3}\right] $-sub-module of 
$\mathcal{F(}E)(3)$ generated by the vectors 
\begin{eqnarray*}
&&\left( x_{1}.\left(
x_{2}.x_{3}\right) -\left( x_{1}.x_{2}\right) .x_{3}\right) ,\,\left(
x_{1}.\left( x_{2}.x_{3}\right) -x_{2}.\left( x_{1}.x_{3}\right) \right)
\,, \\
&&\left( x_{1}.\left( x_{2}.x_{3}\right) -\left( x_{2}.x_{1}\right)
.x_{3}\right) 
\end{eqnarray*}
 for all $x_{1},x_{2},x_{3}\in E$.

Likewise
\begin{eqnarray*}
&&R_{4}^{\perp }= \langle ( x_{1}.( x_{2}.x_{3}) -( x_{1}.x_{2}).x_{3}) ,(
x_{1}.( x_{2}.x_{3}) -x_{3}.(x_{2}.x_{1}) ) , \\ 
&& x_{1}.(x_{2}.x_{3}) -( x_{3}.x_{2}) .x_{1}) \rangle \\ 
&&R_{5}^{\perp }=\langle ( x_{1}.(x_{2}.x_{3}) -( x_{1}.x_{2}).x_{3})
,(x_{1}.( x_{2}.x_{3}) -x_{2}.( x_{3}.x_{1})) , \\ && x_{1}.( x_{2}.x_{3}) -(
x_{2}.x_{3}).x_{1}) ,( x_{1}.( x_{2}.x_{3}) -x_{3}.(x_{1}.x_{2})),( x_{1}.(
x_{2}.x_{3}) -(x_{3}.x_{1}).x_{2}) \rangle  \end{eqnarray*}
and $\dim R_{4}^{\perp }=9, \ \dim R_{5}^{\perp }=10.$ Then it is sufficient to
note that ($\mathcal{F(}E)/\mathcal{R} )^{!}$ is by definition
the binary quadratic operad $\mathcal{F(}E)/\mathcal{R }^{\perp }$.

\subsection{Koszul duality}

Recall that a quadratic operad  $\mathcal{P}$ is called Koszul operad if for every $\mathcal{P}$-free algebra $F_{\mathcal{P}}(V)$ one has 
$H_{i}^{\mathcal{P}}(F_{\mathcal{P}}(V))=0, \quad i > 0$.

\begin{proposition}
The operads   $\mathcal{A}ss$, $\mathcal{V}inb$, $\mathcal{P}reLie$ are Koszul operads. The operads $G_{4}-\mathcal{A}ss$ 
and  $G_{5}-\mathcal{A}ss$ are not Koszul operads.
\end{proposition}
{\it Proof.} In fact, from [G.K] and [G] the operads $\mathcal{A}ss$, $\mathcal{P}reLie$ are Koszul operads. Considering the 
relations between $\mathcal{P}reLie$ and $\mathcal{V}inb$, this last satisfies the same property. 
For the two others, we will show that there are not Koszul operads
using their Poincare series and the Ginzburg-Kapranov criterium. This series is defined for a general operad
 $\mathcal{P}$ by
\[
g _{\mathcal{P}} (x):= \sum_{i=1}^{\infty} (-1)^n dim{\mathcal{P}}(n) \frac{x^n}{n!}, 
\]
The Poincare series for a Koszul operad and for its dual operad are connected by the fonctional equation 
\[
g_{\mathcal{P}}(g_{\mathcal{P} ^!}(x))=x.
\]
But we have 
$$
g_{\mathcal{G}_{4}-\mathcal{A}ss} (x)= -x+x^2-\frac{3}{2}x^3+\frac{59}{4!}x^4+... \quad , \quad g_{{\mathcal{G}_{4}-\mathcal{A}ss}^{!}} (x)= -x+x^2-\frac{1}{2}x^3-\frac{1}{4}x^4+... 
$$
$$
g_{\mathcal{G}_{5}-\mathcal{A}ss} (x)= -x+x^2-\frac{10}{3!}x^3+\frac{39}{4!}x^4+... \quad , \quad g_{{\mathcal{G}_{5}-\mathcal{A}ss}^{!}} (x)= -x+x^2-\frac{1}{3}x^3+\frac{1}{12}x^4+...$$ 
These series do not satisfy the fonctional equation. This proves the proposition.

\section{The categories $G_i$-ASS}
\subsection{Functorial correspondance}
Let $LIE$ and $LIE-AD$ be the categories of Lie algebras and Lie-admissible algebras on $\Bbb{K}$. The correspondance 
$$
T : LIE-AD \longrightarrow LIE
$$ is a functor of which dual functor corresponds to the application
$$
T^{*} : LIE \longrightarrow LIE-AD
$$
with $T^{*}(\mu)=\frac{1}{2} \mu $ (we identifie a Lie algebra with its law).
\subsection{The functor $A \otimes -$}
It is well known that the category $ASS$ of associative algebras on $\Bbb{K}$ is a tensorial category. 
This property is not true for $LIE-AD$ and for all the subcategories $G_i -ASS$ for $i \neq 1$. In this section we will prove that the functor $A \otimes -$ determines a duality between the algebras over the operad $G_{i}-\mathcal{A}ss$ and algebras on the dual operad $G_{i}-\mathcal{A}ss^{! }$.
\begin{theorem}
Let $A$ be a $G_i$-associative algebra. Then  $A \otimes -$ is a covariant functor
$$ A \otimes - : (G_i -ASS)^{!} \longrightarrow G_i -ASS $$
where $ (G_i -ASS)^{!}$  is the category of algebras corresponding to the algebras on the dual operad $G_{i}-\mathcal{A}ss^{!}$. 
\end{theorem}
{\it Proof.} Let $B$ be an algebra. Then $A \otimes B$ is an algebra for the classical product 
$(a_1 \otimes b_1).(a_2 \otimes b_2)=a_1 .a_2 \otimes b_1 .b_2$. This product is in  $G_{i}-\mathcal{A}ss$ 
if and only if $B$ is an algebra on $G_{i}-\mathcal{A}ss^{!}$. In fact
$$((a_1 \otimes b_1).(a_2 \otimes b_2)).(a_3 \otimes b_3)=(a_1.a_2).a_3 \otimes (b_1.b_2).b_3 =(a_1.a_2).a_3 \otimes b_1.b_2.b_3 $$
because $B$ is associative. Let $\mu$ is the law of the algebra $A \otimes B$ and put $X_i=a_i \otimes b_i$.

\noindent First, let us prove the theorem for the category $G_2-ASS$. In this case $A$ is a Vinberg algebra. Then

\noindent $\mu (X_1,\mu (X_2,X_3))-\mu (X_2,\mu (X_1,X_3))-\mu (\mu (X_1,X_2),X_3))+\mu (\mu (X_2,X_1),X_3)) $

\noindent $ =a_1.(a_2.a_3) \otimes b_1.b_2.b_3 - a_2.(a_1.a_3) \otimes b_2.b_1.b_3 -(a_1.a_2).a_3 \otimes b_1.b_2.b_3 $

\noindent $+ (a_2.a_1).a_3 \otimes b_2.b_1.b_3 $

\noindent $= (a_1.(a_2.a_3)- a_2.(a_1.a_3) -(a_1.a_2).a_3 + (a_2.a_1).a_3) \otimes b_1.b_2.b_3 $

\noindent because $B$ is in the category $G_{2}-ASS^!$,

\noindent $ = 0 \otimes b_1.b_2.b_3 $

\noindent because $A$ is a Vinberg algebra

\noindent $ = 0.$ Then $A \otimes B$ also is a Vinberg algebra.

\noindent The demonstration is the same for the other $(G_{i}-\mathcal{A}ss)$-algebras, taking the adapted relation (i.e for $i=3$, we consider
$\mu (\mu (X,Y),Z)-\mu (X,\mu (Y,Z))=\mu (\mu (X,Z),Y)-\mu (X,\mu (Z,Y))$, etc...). 

\medskip 

\noindent {\bf Applications.} This theorem permits to construct interesting classes of Vinberg algebras and to give new 
examples of Lie algebras provided with affine structure. For example suppose $dim A = dim B = 2$. From [R], 
$A$ is isomorphic to one of the following algebras :

\noindent 1. $A$ is commutative and isomorphic to

$$
A_{1}:\left\{ 
\begin{array}{l}
X_{1}.X_{1}=X_{1} \\ 
X_{1}.X_{2}=X_{2}.X_{1}=X_{2} \\ 
X_{2}.X_{2}=X_{2}
\end{array}
\right. \qquad A_{2}:\left\{ 
\begin{array}{l}
X_{1}.X_{1}=X_{1} \\ 
X_{1}.X_{2}=X_{2}.X_{1}=X_{2} \\ 
X_{2}.X_{2}=0
\end{array}
\right. \\ 
$$

$$
A_{3}:\left\{ 
\begin{array}{l}
X_{1}.X_{1}=X_{1} \\ 
X_{1}.X_{2}=X_{2}.X_{1}=X_{2} \\ 
X_{2}.X_{2}=-X_{1}
\end{array}
\right. \\ 
\begin{array}{l}
\qquad A_{4}:\left\{ X_{1}.X_{1}=X_{2}
\right.
\end{array}
$$

$$
\quad \quad
A_{5}:\left\{ 
\begin{array}{l}
X_{1}.X_{1}=X_{1}
\qquad \qquad
\end{array}
\right.
\qquad \quad A_{6}:\left\{ 
\begin{array}{l}
X_{i}.X_{j}=0
\qquad
\quad
\end{array}
\right.
$$

\noindent 2. $A$ is non commutative and isomorphic to

$$
A_{7}:\left\{ 
\begin{array}{l}
X_{1}.X_{1}=\frac{b^{2}+2e}{e}X_{1}-b\frac{b^{2}+e}{e^{2}}X_{2} \\ 
X_{1}.X_{2}=bX_{1}-\frac{b^{2}-e}{e}X_{2} \\ 
X_{2}.X_{1}=bX_{1}-\frac{b^{2}}{e}X_{2} \\ 
X_{2}.X_{2}=eX_{1}-bX_{2}
\end{array}
\right.  \quad A_{8}:\left\{ 
\begin{array}{l}
X_{1}.X_{1}=aX_{1}+cX_{2} \\ 
X_{1}.X_{2}=X_{2} \\ 
X_{2}.X_{1}=X_{2}.X_{2}=0%
\end{array}
\right.
$$
$$
\qquad A_{9}:\left\{ 
\begin{array}{l}
X_{1}.X_{1}=aX_{1} \\ 
X_{1}.X_{2}=(a+1)X_{2} \\ 
X_{2}.X_{1}=aX_{2} \\ 
X_{2}.X_{2}=0
\end{array}
\right.
$$

\noindent In this case the Lie algebra associated to $A$ is the two-dimensional
solvable abelian Lie algebra.

\medskip

\noindent Let us classify the $Vinb^{!}$ algebra of dimension 2.

\noindent 1. $B$ is commutative and isomorphic to $A_{i}$, $i=1,..,7.$

\noindent 2. $B$ is non commutative and isomorphic to
$$
B_{7}:\left\{ 
\begin{array}{l}
e_{1}.e_{1}=e_{1} \\ 
e_{1}.e_{2}=e_{2} \\ 
e_{2}.e_{1}=e_{2}.e_{2}=0
\end{array}
\right.
$$

\noindent If $A$ and $B$ are commutative, the corresponding Lie algebra is the
4-dimensional abelian Lie algebra.

\noindent Suppose $A$ is commutative and $B$ is not commutative. In this case the
bracket of the Lie algebra associated to $A\otimes B$ satisfy%
$$
\left[ X_{i}\otimes e_{j},X_{k}\otimes e_{l}\right] =X_{i}X_{k}\otimes
e_{j}e_{l}-X_{k}X_{i}\otimes e_{l}e_{j}=X_{i}X_{k}\otimes \left[ e_{j},e_{l}%
\right]
$$
with $\left[ e_{1},e_{2}\right] =e_{2}.$

\noindent When we put $f_{ij}=X_{i}\otimes e_{j}$, we obtain the following list of Lie
algebras, $\mathfrak{g}_{i7},$underlying to the Vinberg algebra 
$A_{i}\otimes B_{7}$:

\smallskip

$
\begin{array}{ll}
g_{17} & \quad \left[ f_{11},f_{12}\right] =f_{12},\quad \left[ f_{11},f_{22}
\right] =f_{22},\quad \left[ f_{12},f_{21}\right] =-f_{22},\quad \left[
f_{21},f_{22}\right] =f_{22} \\ 
g_{27} & \quad \left[ f_{11},f_{12}\right] =f_{12},\quad \left[ f_{11},f_{22}
\right] =f_{22},\quad \left[ f_{12},f_{21}\right] =-f_{22} \\ 
g_{37} & \quad \left[ f_{11},f_{12}\right] =f_{12},\quad \left[ f_{11},f_{22}
\right] =f_{22},\quad \left[ f_{12},f_{21}\right] =-f_{22},\quad \left[
f_{21},f_{22}\right] =-f_{12} \\ 
g_{47} & \quad \left[ f_{11},f_{12}\right] =f_{22} \\ 
g_{57} & \quad \left[ f_{11},f_{12}\right] =f_{12} \\ 
g_{67} & \quad abelian
\end{array}
$

\smallskip

\noindent Likewise, if $A$ is a non commutative Vinberg algebra and $B$ a commutative $%
Vinb^{!}$. algebra then the bracket of the corresponding Lie algebra satisfies%
$$
\left[ X_{i}\otimes e_{j},X_{k}\otimes e_{l}\right] =X_{i}X_{k}\otimes
e_{j}e_{l}-X_{k}X_{i}\otimes e_{l}e_{j}=\left[ X_{i},X_{k}\right] \otimes
e_{j}e_{l}
$$
For the algebras $A_{i}\otimes B_{j}$, $i=7,8,9,$ let us note that the
corresponding Lie algebras $\mathfrak{g}_{ij}$ satisfy $\mathfrak{g}_{7j}=%
\mathfrak{g}_{8j}=\mathfrak{g}_{9j}$ for $j=1,...,7.$ Using the same
previous notations we obtain the following Lie algebras :

\smallskip

$
\begin{array}{ll}
g_{i1} & \left[ f_{11},f_{21}\right] =f_{21},\quad \left[ f_{11},f_{22}%
\right] =f_{22},\quad \left[ f_{12},f_{21}\right] =f_{22},\quad \left[
f_{12},f_{22}\right] =f_{22}\quad \quad \\ 
g_{i2} & \left[ f_{11},f_{21}\right] =f_{21},\quad \left[ f_{11},f_{22}%
\right] =f_{22},\quad \left[ f_{12},f_{21}\right] =f_{22} \\ 
g_{i3} & \left[ f_{11},f_{21}\right] =f_{22},\quad \left[ f_{11},f_{22}%
\right] =f_{22},\quad \left[ f_{12},f_{21}\right] =f_{22},\quad \left[
f_{12},f_{22}\right] =-f_{21} \\ 
g_{i4} & \left[ f_{11},f_{21}\right] =f_{22} \\ 
g_{i5} & \left[ f_{11},f_{21}\right] =f_{21} \\ 
g_{i6} & abelian
\end{array}
$

\smallskip

\noindent At last, let us look the case $A=A_{i}$, $i=7,8,9$ and $B$ $=B_{7}.$ We
obtain :

$
\begin{array}{ll}
g_{77}: & \left[ f_{11},f_{12}\right] =\frac{b^{2}+2e}{e}f_{12}-b\frac{%
b^{2}+e}{e^{2}}f_{22},\quad \left[ f_{11},f_{21}\right] =f_{21}, \\ 
& \left[ f_{11},f_{22}\right] =bf_{12}-\frac{b^{2}-e}{e}f_{22}\quad  \\ 
& \left[ f_{12},f_{21}\right] =-bf_{12}+\frac{b^{2}}{e}f_{22},\quad \left[
f_{21},f_{22}\right] =ef_{12}-bf_{22} \\ 
g_{87}: & \left[ f_{11},f_{12}\right] =af_{12},\quad \left[ f_{11},f_{21}%
\right] =f_{21},\quad \left[ f_{11},f_{22}\right] =\left( a+1\right) f_{22},
\\ 
& \left[ f_{12},f_{21}\right] =-af_{22}. \\ 
g_{97}: & \left[ f_{11},f_{12}\right] =af_{12}+cf_{22},\quad \left[
f_{11},f_{21}\right] =f_{21},\quad \left[ f_{11},f_{22}\right] =f_{22}\quad
\quad 
\end{array}
$

\smallskip

\noindent Comparing with the classification of 4-dimensional real Lie algebras
presented in [V], we obtain :

\begin{theorem}
The following solvable Lie algebras have an affine structure of tensorial
type :

$\mathfrak{g}_{3,2}(1)\oplus \mathbb{R},\quad \mathfrak{g}_{4,5}(1,a),\quad 
\mathfrak{g}_{4,6}(1),\quad \mathfrak{g}_{4,9}(\alpha ),\quad \mathfrak{g}%
_{4,10},\quad \mathfrak{g}_{2}\oplus \mathfrak{g}_{2},\quad \mathfrak{g}%
_{4,1}$

$\mathfrak{g}_{4,2},\quad \mathfrak{g}_{2}\oplus \mathbb{R}^{2},\quad 
\mathfrak{g}_{3,1}(1)\oplus \mathbb{R},\quad \mathbb{R}^{4}.$
\end{theorem}

\section{Cohomology of Lie-admissible algebras}

\subsection{Cohomology of ${\mathcal{L}ieAdm}$-algebras}
Let $\mathcal{P}$ be a binary quadratic operad. A $\mathcal{P}$-algebra is given by a $\Bbb{K}$-vector space 
$V$ and an operad morphism
$$ \varphi : \mathcal{P} \longrightarrow \mathcal{E}nd(V) $$
where $\mathcal{E}nd(V)$ is the operad of endomorphismes of $V$ defined by 
$$\mathcal{E}nd(V)(n)=Hom_{\Bbb{K}} (V^{\otimes ^n} , V).$$
Let $A$ be a ${\mathcal{L}ieAdm}$-algebra. The cochain complex of the cohomology 
$H^{*}_{\mathcal{L}ieAdm}(A,A)$ is given by 
$$C^{n}_{\mathcal{L}ieAdm}(A) = Hom_{\Bbb K} ((\mathcal{L}ieAdm ^{!})^{\vee} (n) \otimes _{\sum _{n}} A^{\otimes ^{n}},A)$$
where $V^{\vee}=V^{*} \otimes _{\sum _n} (sgn)$ for $V$ a ${\Bbb K}$-vector space provided with an action of the group $\sum _{n}$, 
$V^{*}=Hom(V,\Bbb{K})$ and $sgn$ the signature representation.

\noindent As every $\mathcal{L}ieAdm ^{!}$-algebra is an associative $3$-order commutative algebra, then the complex 
$C^{*}_{\mathcal{L}ieAdm}(A)$ is 
$$ A \longrightarrow Hom( A \otimes A,A) \longrightarrow Hom(\wedge ^{3}(A),A) \longrightarrow  Hom(\wedge ^{4}(A),A) \longrightarrow ...$$
The differential operator is defined by the composition mapping of the operad $\mathcal{L}ieAdm$. To determinate its expression, let us
recall the definition of Nijenhuis-Gerstenhaber products.

Let $f$ and $g$ be $n$-, respectively $m$-, linear mappings on a vector space 
$V$. We define the $(n+m-1)$-linear mapping $f\odot _{6}g$ by 
\[
\begin{tabular}{l}
$f\odot _{6}g(X_{1},...,X_{n+m-1})=\sum_{i=1,..,n}\sum_{\sigma \in
\sum_{n+m-1}}(-1)^{\varepsilon (\sigma )}(-1)^{(i-1)(m-1)}f(X_{\sigma
(1)},..,$ \\ 
$X_{\sigma (i-1)},g(X_{\sigma (i)},..,X_{\sigma (i+m-1)}),X_{\sigma
(i+m)},..,X_{\sigma (n+m-1)})$%
\end{tabular}
\]
where $\sum_{p}$ refers to the $p-$symmetric group. This product corresponds to the composition product in the operad
$\mathcal{L}ieAdm$.

\noindent For example if $f=g=\mu $ and $n=m=2$ then 
\begin{eqnarray*}
&&\mu \odot_{6}\mu (X_{1},X_{2},X_{3})= \\
&&\sum_{\sigma \in \sum_{3}}(-1)^{\varepsilon (\sigma )}\{\mu (\mu
(X_{\sigma (1)},X_{\sigma (2)}),X_{\sigma (3)})-\mu (X_{\sigma (1)},\mu
(X_{\sigma (2)},X_{\sigma (3)})\}
\end{eqnarray*}
 and $\mu $ is a Lie-admissible law if and
only if $\mu \odot _{6}\mu =0.$ Likewise if $\mu $ is a bilinear
mapping and $f$ and endomorphism of $V$, then 
$$\mu \odot
_{6}f(X_{1},X_{2})=\mu (f(X_{1}),X_{2})-\mu (f(X_{2}),X_{1})+\mu
(X_{1},f(X_{2}))-\mu (X_{2},f(X_{1}))$$
 and $$f\odot _{6}\mu
(X_{1},X_{2})=f(\mu (X_{1},X_{2}))-f(\mu (X_{2},X_{1}))=f([X_{1},X_{2}]_{\mu
})$$ as soon as $\mu \odot _{6}\mu =0.$

In [N], Nijenhuis defines the following product denoted  by $f\overline{o}g$ :
$$
\begin{tabular}{l}
$f\odot
_{1}g(X_{1},...,X_{n+m-1})=%
\sum_{i=1,..,n}(-1)^{(i-1)(m-1)}f(X_{1},..,X_{i-1},g(X_{i},..,X_{i+m-1})$,
\\ 
$X_{i+m},..,X_{n+m-1})$ .
\end{tabular}
$$
We have 
$$f\odot
_{6}g=\sum_{\sigma \in \sum_{n+m-1}} (-1)^{\epsilon (\sigma)}g\odot _{1}f \circ \sigma =\sum_{\sigma \in
\sum_{n+m-1}}(-1)^{\epsilon (\sigma)} f\overline{o}g \circ \sigma .$$ Let $P$ be the antisymmetric operator. It is
defined by 
\[
P(f)(X_{1},...,X_{n})=\sum_{\sigma \in \sum_{n}}(-1)^{\varepsilon (\sigma
)}f(X_{\sigma (1)},X_{\sigma (2)},...,X_{\sigma (n)}). 
\]
It is clear that $f\odot _{6}g=P(f\odot _{1}g)$

\begin{lemma}
We have the following identities : 
\[
P(P(f)\odot _{1}g)=(n+m-1)!P(f\odot _{1}g)=P(f\odot _{1}P(g)). 
\]
\end{lemma}

\noindent This can be prove directly. 

\noindent 
We deduce that the following bracket 
\[
\lbrack f,g]^{\odot _{6}}=f\odot
_{6}g-(-1)^{(n-1)(m-1)}g\odot _{6}f. 
\]
satisfies : 
\[
\begin{tabular}{l}
$1.\quad $ $[g,f]^{\odot _{6}}=(-1)^{(n-1)(m-1)+1}[f,g]^{\odot
_{6}}$ \\ 
$2.\quad $ $(-1)^{(n-1)(p-1)}[[f,g]^{\odot _{6}},h]^{\odot
_{6}}+(-1)^{(m-1)(n-1)}[[g,h]^{\odot _{6}},f]^{\odot _{6}}+$
\\ 
$(-1)^{(p-1)(m-1)}[[h,f]^{\odot _{6}},g]^{\odot _{6}}=0$%
\end{tabular}
\]
where $h$ is a $p$-linear mapping on $V$. 

\smallskip

\noindent {\bf{Consequence}}.
 $[,]^{\odot _6}$ is a bracket of a graded Lie algebra (on the
space of multilinear mappings).

\begin{definition}
Let $\phi$ be a in $C^{n}_{\mathcal{L}ieAdm}(A)$, where $A$ is a $\mathcal{L}ieAdm$-algebra 
of law $\mu$. The differential operator 
$$ \delta _{\mu} : C^{n}_{\mathcal{L}ieAdm}(A) \longrightarrow C^{n+1}_{\mathcal{L}ieAdm}(A)$$
is defined by 
$$ \delta _{\mu} \phi = - [\mu, \phi]^ {\odot _6} $$
\end{definition}

\subsection{Particular cases}

i) $n=2$

\noindent Let $\varphi $ be a bilinear mapping. Then 

$$
\begin{array}{l}

-\frak{\delta }_{\mu }\varphi \left( X_{1},X_{2},X_{3}\right) =   \\
 \mu \left( \varphi \left( X_{1},X_{2}\right) ,X_{3}\right) -\mu \left( X_{1},\varphi
\left( X_{2},X_{3}\right) \right) +\mu \left( \varphi \left(
X_{2},X_{3}\right) ,X_{1}\right)  \\
-\mu \left( X_{2},\varphi \left( X_{3},X_{1}\right) \right) -\mu \left(
X_{3},\varphi \left( X_{1},X_{2}\right) \right) +\mu \left( \varphi \left(
X_{3},X_{1}\right) ,X_{2}\right) \\
-\mu \left( \varphi \left( X_{2},X_{1}\right) ,X_{3}\right) +\mu \left(
X_{2},\varphi \left( X_{1},X_{3}\right) \right) -\mu \left( \varphi \left(
X_{3},X_{2}\right) ,X_{1}\right) \\
+\mu \left( X_{3},\varphi \left( X_{2},X_{1}\right) \right) +\mu \left(
X_{1},\varphi \left( X_{3},X_{2}\right) \right) -\mu \left( \varphi \left(
X_{1},X_{3}\right) ,X_{2}\right) \\
+\varphi \left( \mu \left( X_{1},X_{2}\right) ,X_{3}\right) -\varphi
\left( X_{1},\mu \left( X_{2},X_{3}\right) \right) +\varphi \left( \mu
\left( X_{2},X_{3}\right) ,X_{1}\right) \\
-\varphi \left( X_{2},\mu \left( X_{3},X_{1}\right) \right) -\varphi
\left( X_{3},\mu \left( X_{1},X_{2}\right) \right) +\varphi \left( \mu
\left( X_{3},X_{1}\right) ,X_{2}\right) \\
-\varphi \left( \mu \left( X_{2},X_{1}\right) ,X_{3}\right) +\varphi
\left( X_{2},\mu \left( X_{1},X_{3}\right) \right) -\varphi \left( \mu
\left( X_{3},X_{2}\right) ,X_{1}\right) \\
+\varphi \left( X_{3},\mu \left( X_{2},X_{1}\right) \right) +\varphi
\left( X_{1},\mu \left( X_{3},X_{2}\right) \right) -\varphi \left( \mu
\left( X_{1},X_{3}\right) ,X_{2}\right)
\end{array}
$$

ii) $n=1$

\noindent Let $f$ be in $\mathcal{C}^{1}.$ Then 
\begin{eqnarray*}
\frak{\delta }_{\mu }f\left( X_{1},X_{2}\right) &=&-\mu \left(
f(X_{1}),X_{2}\right) -\mu \left( X_{1},f(X_{2})\right) +f\left( \mu \left(
X_{1},X_{2}\right) \right) \\
&&+\mu \left( f(X_{2}),X_{1}\right) +\mu \left( X_{2},f(X_{1})\right)
-f\left( \mu \left( X_{2},X_{1}\right) \right)
\end{eqnarray*}

\noindent iii) $n=0$

\noindent Here we can define directely the definition of $\delta _{\mu} (X)$.
Let $X$ be in $A$ . Consider the map 
\[
h_{X}:Y\mapsto \mu \left( X,Y\right) -\mu \left( Y,X\right) 
\]
\begin{eqnarray*}
-\frak{\delta }_{\mu }h_{X}\left( X_{1},X_{2}\right) &=&\mu \left(
h_{X}(X_{1}),X_{2}\right) +\mu \left( X_{1},h_{X}(X_{2})\right)-
h_{X}\left( \mu \left( X_{1},X_{2}\right) \right) \\
&&-\mu \left( h_{X}(X_{2}),X_{1}\right)- \mu \left( X_{2},h_{X}(X_{1})\right)
+h_{X}\left( \mu \left( X_{2},X_{1}\right) \right) \\
&=&\mu \left( \mu \left( X,X_{1}\right) ,X_{2}\right) -\mu \left( \mu \left(
X_{1},X\right) ,X_{2}\right) +\mu \left( X_{1},\mu \left( X,X_{2}\right)
\right) \\
&&-\mu \left( X_{1},\mu \left( X_{2},X\right) \right) -h_{X}\left( \mu \left(
X_{1},X_{2}\right) \right) +h_{X}\left( \mu \left( X_{2},X_{1}\right)
\right) \\ &&-\mu \left( \mu \left( X,X_{2}\right) ,X_{1}\right) +\mu \left(
\mu \left( X_{2},X\right) ,X_{1}\right) -\mu \left( X_{2},\mu \left(
X,X_{1}\right) \right) \\
&&+\mu \left( X_{2},\mu \left( X_{1},X\right) \right)
\end{eqnarray*}
\begin{eqnarray*}
-\frak{\delta }_{\mu }h_{X}\left( X_{2},X_{1}\right) &=&\mu \left( \mu \left(
X,X_{2}\right) ,X_{1}\right) -\mu \left( \mu \left( X_{2},X\right)
,X_{1}\right) +\mu \left( X_{2},\mu \left( X,X_{1}\right) \right) \\
&&-\mu \left( X_{2},\mu \left( X_{1},X\right) \right) -h_{X}\left( \mu \left(
X_{2},X_{1}\right) \right) +h_{X}\left( \mu \left( X_{1},X_{2}\right) \right)
\\ &&-\mu \left( \mu \left( X,X_{1}\right) ,X_{2}\right) +\mu \left( \mu \left(
X_{1},X\right) ,X_{2}\right) -\mu \left( X_{1},\mu \left( X,X_{2}\right)
\right) \\
&&+\mu \left( X_{1},\mu \left( X_{2},X\right) \right)
\end{eqnarray*}
Then
\[
\frak{\delta }_{\mu }h_{X}\left( X_{2},X_{1}\right) -\frak{\delta }_{\mu
}h_{X}\left( X_{1},X_{2}\right) =\mu \frak{\odot }_{6}\mu \left(
X,X_{1},X_{2}\right) =0 
\]
\[
\frak{\delta }_{\mu }h_{X}\left( X_{1},X_{2}\right) =\frak{\delta }_{\mu
}h_{X}\left( X_{2},X_{1}\right) 
\]
Let us consider
\[
\mathcal{C}^{0}=\left\{ X\in {\mathcal{A}} \quad / \quad P\left( \delta
_{\mu }h_{X}\right) =\delta _{\mu }h_{X}\right\} 
\]
For $X\in \mathcal{C}^{0},\quad \delta _{\mu }h_{X}=0$. Then we can
define $B^{1}\left( \mathcal{A},\mathcal{A} \right) $ putting 
\[
\delta \left( X\right) =h_{X} 
\]
and $Z^{1}\left( \mathcal{A},\mathcal{A} \right)
=\left\{ f\in \mathcal{C}^{1}/\delta f=0\right\} $

\noindent Then $H^{0}\left( \mathcal{A},\mathcal{A}\right) $ is well defined.

\medskip

\noindent {\bf Remark.} In the way, we can define the cohomology $H^{*}_{G_{i}-\mathcal{A}ss}(A,A)$ for a $G_{i}\mathcal{A}ss$-algebra. We denote by
$f\odot_{i}g$ the corresponding Nijenhuis product. We have yet given the expression of this product for $i=1$ and $i=6$. In other cases we put :
$$
\begin{tabular}{l}
$f\odot _{2}g(X_{1},...,X_{n+m-1})=\sum_{\sigma \in
\sum_{n+m-2}}(-1)^{\varepsilon (\sigma
)}\{\sum_{i=1,..,n-1}(-1)^{(i-1)(m-1)}f(X_{\sigma (1)},$ \\ 
$..,$ $X_{\sigma (i-1)},g(X_{\sigma (i)},..,X_{\sigma (m+i-1)}),X_{\sigma
(m+i)},.,X_{\sigma (n+m-2)},X_{n+m-1})$ $+$ \\ 
$(-1)^{(n-1)}f(X_{\sigma (1)},..,X_{\sigma (n-1)},g(X_{\sigma
(n)},..,X_{\sigma (m+n-2)},X_{n+m-1})\}.$%
\end{tabular}
$$

\[
\begin{tabular}{l}
$f\odot _{3}g(X_{1},...,X_{n+m-1})=\sum_{\sigma \in
\sum_{n+m-2}}(-1)^{\varepsilon (\sigma )}\{f(g(X_{1},X_{\sigma
(2)},..,X_{\sigma (m)}),$ \\ 
$X_{\sigma (m+1)},..,X_{\sigma
(m+n-1)})+\sum_{i=2,..,n-1}(-1)^{(i-1)(m-1)}f(X_{1},$ $X_{\sigma (2)},..,$ $%
X_{\sigma (i)},$ \\ 
$g(X_{\sigma (i+1)},..,X_{\sigma (m+i)}),.,X_{\sigma (n+m-1)})\}$%
\end{tabular}
\]

\[
\begin{tabular}{l}
$f\odot _{4}g(X_{1},...,X_{n+m-1})=\sum_{\sigma \in
\sum_{n+m-2}}\sum_{i=0,..,n-1}(-1)^{\varepsilon (\sigma )+i(m-1)}f(X_{\sigma
(1)},X_{2},$ \\ 
$X_{\sigma (3)},...,X_{\sigma (i)},g(X_{\sigma (i+1)},..,X_{\sigma
(m+i)}),...,X_{\sigma (n+m-1)}).$%
\end{tabular}
\]

\[
\begin{tabular}{l}
$f\odot _{5}g(X_{1},...,X_{n+m-1})=\sum_{\sigma \in
A_{n+m-1}}\sum_{i=0,..,n-1}(-1)^{i(m-1)}f(X_{\sigma (1)},X_{\sigma (2)},$ \\ 
$X_{\sigma (3)},...,X_{\sigma (i)},g(X_{\sigma (i+1)},..,X_{\sigma
(m+i)}),...,X_{\sigma (n+m-1)}).$%
\end{tabular}
\]
where $A_{p}$ refers to the alternated group.

\subsection{Lie-admissible cohomology of Lie algebras}
If $\frak{g}$ is a Lie algebra of law $\mu$, it is also a Lie-admissible algebra. Then it is possible to consider the following cohomologies

1. $H^{*}_{\mathcal{L}ieAdm} (\frak{g},\frak{g})$ 

2. $H^{*}_{\mathcal{L}ie} (\frak{g},\frak{g)}$ 

\noindent in these two cases, the cochain spaces are the same, excepted for $n=2$. Recall also that the second cohomoly in nothing other that the Chevalley cohomology.

\begin{definition}
Let $\frak{g}$ be a Lie algebra. We call Lie-admissible cohomology of $\frak{g}$ with value in $\frak{g}$ the cohomology $H^{*}_{\mathcal{L}ieAdm}(\frak{g},\frak{g})$.
\end{definition}

\noindent {\bf Remark.}

\noindent  A 2-cochain
corresponding to the Chevalley's cohomology is alternated and satisfies
\begin{eqnarray*}
\delta ^{c}\varphi \left( X_{1},X_{2},X_{3}\right) &=&\varphi \left(
X_{1},X_{2}\right) .X_{3}+\varphi \left( X_{2},X_{3}\right) .X_{1}+\varphi
\left( X_{3},X_{1}\right) .X_{2} \\
&&+\varphi \left( X_{1}.X_{2},X_{3}\right) +\varphi \left(
X_{2}.X_{3},X_{1}\right) +\varphi \left( X_{3}.X_{1},X_{2}\right).
\end{eqnarray*}
where $\delta ^{c}$ is the Chevalley operator.
If we compute $\delta \varphi \in B^{2}_{\mathcal{L}ieAdm}(\frak{g},\frak{g})$ we obtain 
$\varphi \left( X,Y\right) =-\varphi \left( Y,X\right) $ then 
\begin{eqnarray*}
\delta \varphi \left( X_{1},X_{2},X_{3}\right) &=&2 (
\varphi \left( X_{1},X_{2}\right) .X_{3}-X_{3}.\varphi \left(
X_{1},X_{2}\right) +\varphi \left( X_{2},X_{3}\right) .X_{1} \\
&&-X_{1}.\varphi\left( X_{2},X_{3}\right)  
+\varphi \left( X_{3},X_{1}\right) .X_{2}-X_{2}.\varphi \left(
X_{3},X_{1}\right) \\
&&+\varphi \left( X_{1}.X_{2},X_{3}\right) -\varphi \left(
X_{3},X_{1}.X_{2}\right)  
+\varphi \left( X_{2}.X_{3},X_{1}\right) \\
&&-\varphi \left(X_{1},X_{2}.X_{3}\right) +\varphi \left(
X_{3}.X_{1},X_{2}\right) -\varphi \left( X_{2},X_{3}.X_{1}\right) 
) 
\end{eqnarray*}
then
\[
\delta \varphi  =4\delta ^{c}\varphi.
 \]

\section{Lie-admissible modules on a Lie algebra}

\subsection{Module on a Lie-admissible algebra}

Let $\mathcal{A=}(A,\mu )$ be a Lie-admissible algebra and $M$ a vector
space on $\Bbb{K}$.

\begin{definition}
$M$ is an $\mathcal{A}$-module if there is bilinear mapping 
\begin{eqnarray*}
\lambda &:&A\otimes M\rightarrow M \\
\rho &:&M\otimes A\rightarrow M
\end{eqnarray*}
satisfying :

$\lambda (X,\lambda (Y,v)-\lambda (Y,\lambda (X,v)-\lambda ([X,Y]_{\mu
},v)-\lambda (X,\rho (v,Y)+\rho (\lambda (X,v),Y)+\rho (v,[X,Y]_{\mu })-\rho
(\rho (v,X),Y)-\rho (\lambda (Y,v),X)+\lambda (Y,\rho (v,X))+\rho (\rho
(v,Y),X)=0$

for all $X,Y\in A$ and $v\in M.$
\end{definition}

For example, the vector space $A$ is an $\mathcal{A}$-module.
\begin{proposition}
Let $\mathcal{A=}(A,\mu )$ be a Lie-admissible algebra and $M$ an $\mathcal{A%
}$-module defined by the mapping $\lambda $ and $\rho .$ Then the bilinear
mapping 
\[
\widehat{\lambda }:A\otimes M\rightarrow M \]
defined by 
\[
\widehat{\lambda }(X,v)=\lambda (X,v)-\rho (v,X) 
\]
provides the vector space $M$ with an $\mathcal{A}_{L}$-module where $%
\mathcal{A}_{L}$ is the Lie algebra $(A,[,]_{\mu }).$
\end{proposition}

We find again the same result established by Nijenhuis in [N] for the
Vinberg algebra.

\subsection{Modules on $G_i$-associative algebras}

If $\mathcal{A}_{V}\mathcal{=}(A,\mu _{V})$ is a Vinberg algebra an $%
\mathcal{A}_{V}$-module $M$ is given by the mappings $\lambda $ and $\rho $
satisfying the two conditions 
\[
\left\{ 
\begin{array}{l}
\lambda (X,\lambda (Y,v)-\lambda (Y,\lambda (X,v)-\lambda ([X,Y]_{\mu },v)=0
\\ 
\lambda (X,\rho (v,Y)-\rho (\lambda (X,v),Y)-\rho (v,\mu (X,Y))+\rho (\rho
(v,X),Y)=0.
\end{array}
\right. 
\]
Considering $\mathcal{A}_{V}$ as a Lie-admissible algebra, then $M$ is also
a module on this Lie-admissible algebra.

The notion of $\mathcal{A}$-module is well known in the cases of 
Lie-admissible algebras of
type 1 (associative). It is easy to
write the definitions of modules on algebras of type 4 and 5. We find
:

- type 4 : 
\[
\left\{ 
\begin{array}{c}
\lambda (\mu (X,Y),v)-\lambda (X,\lambda (Y,v)-\rho (\rho (v,Y),X)+\rho
(v,\mu (Y,X))=0 \\ 
\rho (\lambda (X,v),Y)-\lambda (X,\rho (v,Y)-\rho (\lambda (Y,v),X)+\lambda
(Y,\rho (v,X))=0
\end{array}
\right. 
\]

-type 5 :

$\lambda (\mu (X,Y),v)-\lambda (X,\lambda (Y,v)+\rho (\lambda
(X,v),Y)-\lambda (X,\rho (v,Y)-\rho (v,\mu (X,Y))+\rho (\rho (v,X),Y)=0.$ 

\noindent We
can see that if $\mu $ is antisymmetric (i.e. a law of Lie algebra) then we
find again the definition of module on Lie algebra considering $\rho
(v,X)=-\lambda (X,v).$

\subsection{Lie-admissible modules on Lie algebras}

Let $\mathcal{A=}(A,\mu )$ be a Lie algebra. Consider this algebra as a Lie-admissible algebra that we will note $\mathcal{A}_{ad}\mathcal{=}(A,\mu )$
to distinguish the two structures. It is clear that every module $M$ on the
Lie algebra $\mathcal{A}$ is also a module on the Lie-admissible algebra $%
\mathcal{A}_{ad}.$ But the converse is false.

\begin{definition}
We call Lie-admissible module on the Lie algebra $\mathcal{A=}(A,\mu )$ every
module on the Lie-admissible algebra $\mathcal{A}_{ad}\mathcal{=}(A,\mu ).$
\end{definition}

Let $\frak{g}$ be the solvable non abelian 2-dimensional Lie algebra. There
exists a basis $\{X_{1},X_{2}\}$ such that $[X_{1},X_{2}]=X_{2}.$ Every one
dimensional module on the Lie algebra $\frak{g}$ is given by the mapping $%
\lambda $ (here $\rho =-\lambda )$ defined by 
\begin{eqnarray*}
\lambda (X_{1},v) &=&\alpha v \\
\lambda (X_{2},v) &=&0
\end{eqnarray*}
On the other hand, a Lie-admissible module on $\frak{g}$ is determined by the
mapping $\lambda $ and $\rho $ given by 
\[
\left\{ 
\begin{array}{c}
\lambda (X_{1},v)=\alpha v \\ 
\lambda (X_{2},v)=\beta v
\end{array}
\right. ;\quad \left\{ 
\begin{array}{c}
\rho (v,X_{1})=\gamma v \\ 
\rho (v,X_{2})=\beta v
\end{array}
\right. 
\]
Suppose now that $M$ is a $n$-dimensional Lie-admissible module on $\frak{g}$%
.\thinspace Then if $A,B,C,D$ are the matrices on the linear operators $%
\lambda (X_{1},.),\lambda (X_{2},.),\rho (.,X_{1}),\rho (.,X_{2})$ in a
given basis of $M$, then these matrices satisfy 
\[
\lbrack (B-D),(C-A)]=B-D 
\]
We describe in this way all the structures of Lie-admissible modules on $\frak{g}$.

Now consider the Lie algebra $\frak{g}=sl(2,\Bbb{C)}$. By a similar
computation we can see that every $n$-dimensional Lie-admissible module on $sl(2,%
\Bbb{C)}$ is described by the following matrix representations : 
\[
\left\{ 
\begin{array}{l}
\lbrack A_{1}-B_{1},A_{2}-B_{2}]=4(A_{2}-B_{2}) \\ 
\lbrack A_{1}-B_{1},A_{3}-B_{3}]=-4(A_{3}-B_{3}) \\ 
\lbrack A_{2}-B_{2},A_{3}-B_{3}]=2(A_{1}-B_{1})
\end{array}
\right. . 
\]
Such representation also is completely reducible.

\section{Deformations} 

Let us denote $\mathcal{L}A_n$ the algebraic variety of $n$-dimensional Lie-admissible algebras on an algebraic closed field  $\Bbb{K}$ 
of characteristic $0$. We have a natural fibration
\[
\pi :\mathcal{L}A_n \rightarrow \mathcal{L}_n,
\]
where $\mathcal{L}_n$ indicates the algebraic variety of $n$-dimensional Lie algebras :
$$
\pi (\mu) = [,]_{\mu}
$$
This fiber owns a global section
$$s:\mathcal{L}_n \rightarrow \mathcal{L}A_n$$
defined by
$$s(\mu)=\frac{1}{2}\mu.$$ 
Let us denote  $T_{s(\mu)}\mathcal{L}A_n$  and $T_{s(\mu)}\pi ^{-1}((\mu))$ the tangent spaces to the variety and to the fiber at the point
$s(\mu)$. We know that $T_{s(\mu)}\mathcal{L}A_n$ is identified to the space of the 2-cocycles $Z^{2}_{\mathcal{L}ieAdm}(A,A)$.

\begin{lemma}
$T_{s(\mu)}\pi^{-1}(\mu)\simeq $
$\left\{ \varphi:A^{\otimes 2} \rightarrow A \,/ \,\varphi(x,y)=\varphi(y,x)\right\} \simeq S^2(A)$.
\end{lemma}
\noindent {\it Proof.} In fact $\frac{1}{2}\mu+t \varphi$ is a linear deformation of $\frac{1}{2}\mu$ in the fiber $\pi ^{-1} (\mu )$ then 
$\varphi(x,y)-\varphi(y,x)=0$. More every symmetric bilinear mapping is a cocycle that is in $Z^{2}_{\mathcal{L}ieAdm}(A,A)$.

\medskip

Recall the geometrical problem concerning the existence of affine structure on solvable Lie algebras. A Lie algebra is provided with an affine structure if and only if
there exists a Vinberg algebra whose associated Lie algebra corresponds to the given. We know that there exists nilpotent Lie algebras without affine structure. In this case the fiber
$\pi ^{-1} (\mu )$  does not cut the subvariety of Vinberg laws.
For $\varphi \in S^2(A) \subset Z_{\mathcal{L}ieAdm}^2(A,A)$ fixed, the staightline $\frac{1}{2}\mu+t \varphi$ is in the fiber 
$\pi ^{-1} (\mu )$.
This line cut the subvariety of Vinberg if and only if there is $t_0$ such that $\frac{1}{2}\mu+t_0 \varphi$ is a Vinberg law. Considering 
$\varphi$ for $t_0 \varphi$ (we can always suppose that $t_0 = 1$ we obtain :

\begin{proposition} 
The deformation $1/2 \mu +\varphi $ is in $Vinb_{n}$ if and only if the 
symmetric mapping satisfies  
\begin{eqnarray*}
&&4\varphi \left( \mu \left( X_{2},X_{1}\right) ,X_{3}\right) +2\mu \left(
X_{1},\varphi \left( X_{2},X_{3}\right) \right)  +2\varphi \left( X_{1},\mu
\left( X_{2},X_{3}\right) \right) \\ 
&& +4\varphi
\left( X_{1},\varphi \left( X_{2},X_{3}\right) \right)  
-2\mu \left( X_{2},\varphi \left( X_{1},X_{3}\right) \right) -2\varphi
\left( X_{2},\mu \left( X_{1},X_{3}\right) \right)  \\
&&-4\varphi \left( X_{2},\varphi \left( X_{1},X_{3}\right) \right) + \mu
\left( \mu \left( X_{2},X_{1}\right) ,X_{3}\right) = 0.
\end{eqnarray*}
for all $X_{1},X_{2},X_{3}\in \frak{g}$.
\end{proposition}

This proposition can be considered as a criterium of existence of an affine structure on a given Lie algebra.
 
\bigskip

{\Large BIBLIOGRAPHY }

\medskip

\noindent [A] Albert A.A., {\it On the power-associative rings}. Trans. Amer.
Math. Soc. {\bf 64}, (1948). 552--593.

\noindent [Ba] Balavoine D., {\it Deformations of algebras over a quadratic
operad. Operads:} Proceedings of Renaissance Conferences (Hartford,
CT/Luminy, 1995), Contemp. Math., {\bf 202} , Amer. Math. Soc., Providence,
RI, (1997),207--234. 

\noindent [Be] Besnoit Y. {\it Une nilvari\'{e}t\'{e} non affine.} 
J.Diff.Geom. {\bf{41}} (1995), 21-52

\noindent [Bu] Burde D. {\it Affine structures on nilmanifolds}. 
Int. J. of Math,  {\bf 7} (1996), 599-616

\noindent [C.L] Chapoton F., Livernet M., {\it Pre-Lie algebra and the rooted
trees operad}. Preprint QA/0002069

\noindent [D.O] Dekimpe K., Ongenae V., {\it On the number of abelian
left symmetric algebras. }P.A.M.S.,  {\bf 128 (11),} 3191-3200, 2000

\noindent [Dz] Dzhumadil'daev A., {\it Cohomologies and deformations of right
symmetric algebras}. J. Math. Sci. {\bf 93}, (1999), 836-876

\noindent[5] Gerstenhaber M., The cohomology structure of an associative ring,
Ann of math.(2) 78 (1963) 267-288.

\noindent [G.K] Ginzburg V., Kapranov M., Koszul duality for operads. {\it
Duke Math Journal} . {\bf 76,1} (1994), 203-272 

\noindent [Go] Gonzalez Santos, Elduque A. {\it Flexible Lie-admissible
algebras with $A$ and $A\sp{-}$ having the same lattice of subalgebras.} 
Algebras Groups Geom. {\bf 1} (1984), 137--143. 

\noindent [L] Loday J.L., La renaissance des op\'erades. (French) [The rebirth
of operads] S\'eminaire Bourbaki, Vol. 1994/95. {\it Ast\'erisque} {\bf 237},
(1996), 47-74. 

\noindent [M] May J.P. {\it Operadic tensor products and smash products.} 
Operads: Proceedings of Renaissance Conferences (Hartford, CT/Luminy, 1995),
287--303, Contemp. Math., 202, Amer. Math. Soc., Providence, RI, 1997. 

\noindent [N] Nijenhuis A., {\it Sur une classe de propri\'et\'es communes
\`{a} quelques types diff\'erents d'alg\`ebres.} Enseignement Math. (2) {\bf
14}, (1970), 225--277 .

\end{document}